\documentclass{gtmon_a}
\pdfoutput=1

\proceedingstitle{Heegaard splittings of 3-manifolds (Haifa 2005)}
\conferencestart{10 July 2005}
\conferenceend{19 July 2005}
\conferencename{Heegaard splittings of 3-manifolds}
\conferencelocation{Haifa}

\editor{Cameron Gordon}
\givenname{Cameron}
\surname{Gordon}

\editor{Yoav Moriah}
\givenname{Yoav}
\surname{Moriah}

\title{Problems around 3--manifolds}

\author{J\,H Rubinstein}  
\givenname{J Hyam}
\surname{Rubinstein}
\address{Department of Mathematics and Statistics\\
University of Melbourne\\\newline
Parkville, 3010\\Australia}
\email{rubin@ms.unimelb.edu.au}
\urladdr{www.ms.unimelb.edu.au/~rubin}


\arxivreference{}  

\volumenumber{12}
\issuenumber{}
\publicationyear{2007}
\papernumber{11}
\startpage{285}
\endpage{298}
\doi{}
\MR{}
\Zbl{}
\received{24 April 2006}
\revised{29 July 2007}
\accepted{29 July 2007}
\published{3 December 2007}
\publishedonline{3 December 2007}
\proposed{}
\seconded{}
\corresponding{}
\version{}

\newtheorem{prob}{Problem}   
\theoremstyle{remark}
\newtheorem{rem}{Remark}

\begin{document}

\begin{abstract}    

This is a personal view of some problems on minimal surfaces, Ricci
flow, polyhedral geometric structures, Haken 4--manifolds, contact
structures and Heegaard splittings, singular incompressible surfaces
after the Hamilton--Perelman revolution.

\end{abstract}

\maketitle


We give sets of problems based on the following themes;

\begin{itemize}

\item Minimal surfaces and hyperbolic geometry of 3--manifolds. In particular, how do minimal surfaces give information about the geometry of hyperbolic 3--manifolds and conversely how does the geometry affect the types of minimal surfaces which can occur?

\item Ricci flow. Here it would be good to be able to visualize Ricci flow and understand more about where and how singularities form. Also this leads into the next topic, via the possibility of a combinatorial version of Ricci flow. 

\item Combinatorial geometric structures. Various proposals have been made for combinatorial versions of differential geometry in dimension 3. However having a good model of sectional or Ricci curvature is a difficult challenge.

\item Haken 4--manifolds. These are a class of 4--manifolds which admit hierarchies, in the same sense as Haken 3--manifolds. More information is needed to determine how abundant they are.

\item Contact structures and Heegaard splittings. Contact structures are closely related to open book decompositions, which can be viewed as special classes of Heegaard splittings. 

\item Singular incompressible surfaces. A well-known and important conjecture is that all (closed or complete finite volume) hyperbolic 3--manifolds admit immersed or embedded closed incompressible surfaces. We give some variants of this problem. 

\end{itemize}

I would like to thank the referee for a number of useful comments, which hopefully have improved the exposition. 

\section{Minimal surfaces and hyperbolic geometry}

Throughout, hyperbolic 3--manifolds will be either closed or complete with finite volume, unless
indicated otherwise. 

\bigskip

\begin{prob}  

{\rm 
A quasi-Fuchsian subgroup of the fundamental group of a hyperbolic 3--manifold can be represented by a minimal (least area) surface by Schoen and Yau \cite{SY}. There are a number of natural questions about such minimal surfaces.} 

\begin{itemize}

\item When is there a unique such surface representing a given quasi-Fuchsian subgroup, or equivalently in its homotopy class?
\end{itemize}

\begin{rem}\rm 
Note that Ben Andrews has an unpublished flow which shows that if
there is any surface in this homotopy class with all principal
curvatures having absolute value strictly less than $1$, then the flow
produces a minimal surface with the same bounds on principal
curvatures. (In fact, the maximum absolute value of principal
curvature for this minimal surface is the smallest possible amongst
all surfaces in the homotopy class). It then follows easily by a
maximum principle argument due to Uhlenbeck \cite{U} that this surface
is indeed unique in its homotopy class. See also Taubes \cite{Ta} and
Krasnov and Schlenker \cite{KS} for more recent results on properties
of minimal surfaces in hyperbolic 3--manifolds, with the latter
reference concentrating on the quasi-Fuchsian case. See also
Rubinstein \cite{Ru1} for some other properties of minimal surfaces,
including construction of examples of parallel minimal surfaces in
hyperbolic 3--manifolds, by using very short closed geodesics as
barriers.

\end{rem}

\begin{itemize}

\item How does the geometry of such minimal surfaces correspond to the
geometry of the corresponding quasi-Fuchsian 3--manifold, in the case
that there is a unique minimal surface in the quasi-Fuchsian homotopy
class? For example, how does the second fundamental form and area of
such a minimal surface correspond to the volume of the convex core of
the quasi-Fuchsian 3--manifold? Could there be a sequence of
quasi-Fuchsian surface groups of fixed genus, with unique minimal
surface and volume of the convex core approaching infinity?

\item The generalised Gauss map (cf Krasnov and Schlenker \cite{KS}
and Epstein \cite{Ep}) for a quasi-Fuchsian 3--manifold, with a minimal
surface with all principal curvatures having absolute value strictly
less than one, defines a special homeomorphism between the two Riemann
surfaces defined by the domains of discontinuity. Is there a way of
using this homeomorphism to relate the distance between the two
Riemann surfaces, using any of the usual distances on Teichm\"uller
space, with properties of the minimal surface?

\end{itemize}

\end{prob}

\bigskip

\begin{prob}[S\,T Yau]

Does every closed Riemannian (or in particular hyperbolic) 3--manifold admit infinitely many closed embedded minimal surfaces?

\begin{rem}\rm  
One could also ask an easier question allowing merely immersions for
the minimal surfaces---cf Rubinstein \cite{Ru1} for some discussion of
this weaker question. Also note the analogy with the closed geodesic
problem in Riemannian geometry, which asks to find infinitely many
simple closed embedded or immersed geodesic loops in a hyperbolic
3--manifold. Note that it is not known even that a hyperbolic
3--manifold always has infinitely many such embedded loops---cf
Kuhlmann \cite{Ku} for some recent work on this latter problem. The techniques
of Kapouleas \cite{Ka} may be relevant to this problem, since if one
could find two embedded intersecting minimal surfaces and perform a
`desingularisation' by adding arbitrarily many handles along the
curves of intersection, then this would give a solution. Note also
that it seems possible that any Heegaard splitting of sufficiently
large genus (not necessarily irreducible or strongly irreducible)
might be realisable by a minimal surface, but this seems to be a very
difficult question.
\end{rem}

\end{prob}

\begin{prob}

Can a hyperbolic 3--manifold, which is a closed surface bundle over a
circle, have the property that every leaf of the bundle structure is a
minimal surface?

 \begin{rem}\rm  Sullivan \cite{Su} constructed metrics which have this property, but they are unlikely to be hyperbolic. Note that if such a hyperbolic 3--manifold $M$ has some very short geodesic loops, as can happen by doing large longitudinal Dehn surgery along a loop which lies in a leaf, then it is straightforward to prove that not all the leaves of $M$ can be minimal. For the barrier argument in \cite{Ru1} together with a simple maximum principle argument gives a contradiction, as minimal leaves cannot touch such a barrier on the convex side. 
 
\end{rem}

\end{prob}

\begin{prob}
Can minimal surfaces be constructed as suitable limits of normal surfaces?

 \begin{rem}\rm  Note if a 3--manifold is embedded isometrically into a high dimensional Euclidean space, by Nash's theorem, then one can take a suitable net of closely spaced points so that the convex hull of subsets of 4 points gives a Euclidean triangulation, which is a good approximation to the Riemannian 3--manifold. For a fine enough net, any particular smooth surface will be well approximated by a normal surface, relative to such a triangulation. Does this give a useful way of producing sequences of triangulations and normal surfaces, which converge to all minimal surfaces?

\end{rem}

\end{prob}

\begin{prob}[H Rosenberg]

R Bryant \cite{Br} showed that surfaces of constant mean curvature $+1$ have beautiful closed form representations in hyperbolic space ${\bf H}^3$. (This is analogous to the Weierstrass representation of minimal surfaces in ${\bf R}^3$). There is a rich theory of complete Bryant surfaces, which is the usual term for surfaces of constant mean curvature $+1$, in hyperbolic geometry.  Show that every closed or complete finite volume hyperbolic 3--manifold has an embedded closed Bryant surface.

 \begin{rem}\rm One approach would be to study a surface of smallest
 area, which divides the 3--manifold up into two regions of fixed
 volume. Such a surface will have constant mean curvature, but there
 are difficulties in establishing the existence of smooth
 solutions. For the 2--sphere, there are several ways of producing a
 simple closed geodesic loop. Two important ones are the minimax
 method constructing a curve of least length dividing the surface into
 two pieces of equal area. See Hass and Scott \cite{HS} and
 Croke \cite{Cr}. One could start by studying the analogous problem of
 finding loops of constant curvature on a Riemannian 2--sphere.

\end{rem}

\end{prob}

\begin{prob}

Find a `topological' description of embedded closed minimal surfaces in hyperbolic 
3--manifolds. 

 \begin{rem}\rm In the seven Thurston's geometries which are non
 hyperbolic, one has a complete understanding of how minimal surfaces
 can occur topologically---essentially as either incompressible
 surfaces or as Heegaard surfaces or suitable generalisations. See
 \cite{Ru1} for a discussion of this. It is known that strongly
 irreducible splittings can be isotoped to be minimal, or to double
 cover a one-sided surface after a single handle pinch. Also
 incompressible surfaces can be isotoped to be minimal. On the other
 hand, by the barrier argument in \cite{Ru1}, it is clear that there
 are many other topological types of minimal surfaces.  Namazi
 \cite{Na} and Namazi and Souto \cite{NS} construct long cylindrical
 regions in hyperbolic 3--manifolds. Do these regions contain many
 disjoint parallel minimal surfaces?
\end{rem}

\end{prob}

 \section{Ricci flow}
 
\begin{prob}

Is there a `weak' solution to Ricci flow which continues through singularities?

 \begin{rem}\rm Brakke \cite{Bk} achieved such a setup for mean
 curvature flow and this has proved to be very useful. For a more
 recent account of generalized mean curvature flow, see Ecker
 \cite{Ec}. In a brief numerical study (Rubinstein and Sinclair
 \cite{RS}), it is observed for rotationally symmetric neck pinching
 that one might be able to continue by allowing indefinite metrics in
 the part of the manifold where the pinching has occurred. If weak
 solutions could be found, then the advantage would be to continue
 flowing on manifolds, without the need to perform surgery. For an
 analysis of singularity formation in the rotationally symmetric case,
 see Simon \cite{Si1} and Angenent and Knopf  \cite{AK}.
 
\end{rem}

\end{prob}
 
\begin{prob}

For an atoroidal Haken 3--manifold, singularities will form away from
a hierarchy consisting of least area incompressible surfaces. Is there
a way of continuing the flow in a neighbourhood of the hierarchy after
such singularities have formed, without surgery?
 
 \begin{rem}\rm Note that Hamilton \cite{Ha} has analysed neck pinching and it is clear that least area incompressible surfaces will stay away from such necks. Moreover Gabai \cite{Ga1} has shown that even though minimal surfaces do not move continuously as metrics are deformed, one can interpolate between their positions in a useful way. So one could study the motion of the hierarchy right up to the moment of a singularity and possibly use this as a guide to continue the flow past the singularity, as all the topological structure of the manifold is contained in the hierarchy. 
 
 \end{rem}

\end{prob}
 
\begin{prob}

 Even if the initial curvature is strictly negative, it is not known if the Ricci flow develops singularities. The cross curvature flow of Chow--Hamilton \cite{CH} has been proposed to show that any negatively curved metric flows directly to a hyperbolic metric. The cross-curvature flow is a special flow which is only defined for 3--manifolds, since by the construction of Gromov--Thurston \cite{GT}, it is known that there are arbitrarily pinched negatively curved metrics on higher dimensional manifolds, which cannot be deformed to be hyperbolic. If the cross curvature flow could be shown to work as above, this would show that the path components of the space of negatively curved metrics on a 
3--manifold are contractible, following Gabai's solution of the Smale conjecture for hyperbolic 3--manifolds \cite{Ga1}.

\end{prob}
 
\begin{prob}
 
Find a combinatorial Ricci flow. To do this, we first need a robust combinatorial definition of sectional curvature or Ricci curvature in dimension~3. This is very challenging. Similarly, is there a combinatorial analogue of Perelman's entropy? Chow and Luo \cite{CL} have studied combinatorial Ricci flow for surfaces. Various proposals have been made for deforming shapes of tetrahedra by flows, but to obtain success, one would need to choose the `right' triangulation to start with. Otherwise a feature of the flow should involve changing the cell decomposition. 
 
\end{prob}
 
\begin{prob}
 
If Ricci flow did not develop singularities, one could deduce a proof
of the Smale conjecture for all geometric 3--manifolds simultaneously!
Is there any way of using Ricci flow with surgery to deduce
information about the homotopy type of the space of diffeomorphisms?
 
 \begin{rem} \rm
 
 Note that the space of Riemannian metrics on a manifold is homotopy
 equivalent to the manifold. Hence if a 3--manifold admits a geometric
 structure, then the space $\cal D$ of diffeomorphisms of the manifold
 to itself induce the collection of all geometric structures $\cal
 D/I$, where $\cal I$ is the group of isometries of the geometric
 structure. If one could retract the metrics onto the geometric
 structures by Ricci flow, this would prove that the diffeomorphism
 group is homotopy equivalent to the isometry group. Various cases of
 the Smale conjecture have been established for classes of geometric
 3--manifolds by Hatcher \cite{Ha1,Ha2,Ha3}, Ivanov \cite{Iv2,Iv1},
 Gabai \cite{Ga1}, McCullough and Rubinstein \cite{McR},
 and Hong, McCullough and Rubinstein \cite{HMcR}.  \end{rem}

\end{prob}
 
\begin{prob}
 
Can we visualise Ricci flow, eg, is there a natural way of embedding
Ricci flow as a submanifold flow in ${\bf R}^N$ for $N$ large? If so
then level set methods (cf Sethian \cite{Se}) could be used for
numerical study of Ricci flow.
 
 \end{prob}
 
 \section {Combinatorial geometric structures}
 
\begin{prob}
  
Find hyperbolic structures on 3--manifolds by solving the gluing equations for a suitable triangulation. In the case of an irreducible atoroidal manifold with tori boundary, the natural choice is an ideal triangulation. 
  
 \begin{rem}\rm
 
  Casson has suggested a flow along the gradient of volume in the
  ideal case, to change the shapes of tetrahedra in a suitable
  way. Feng Luo \cite{Lu} has shown that a similar flow can be defined
  in the closed case. Lackenby \cite{La} used taut foliations and
  sutured hierarchies (Gabai \cite{Ga}, Scharlemann \cite{Sc}) to
  construct taut ideal triangulations in the ideal case.  Kang and
  Rubinstein \cite{KR} showed that these deform to angle structures if
  and only if a certain normal surface theory obstruction vanishes for
  the triangulation. The challenge is to find ways of modifying a taut
  triangulation to eliminate this obstruction, to at least always
  obtain an angle structure.  In the closed case, a related question
  is whether any hyperbolic 3--manifold admits a one vertex straight
  triangulation, i.e. all edges would then be geodesic loops with a
  possible corner at the single vertex? Properties of these
  triangulations are also interesting---for example are these
  1--efficient (cf Jaco and Rubinstein \cite{JR1,JR2})?  \end{rem}

\end{prob}

\begin{prob}

If an atoroidal 3--manifold has a cubing of non-positive curvature, does Ricci flow immediately produce a metric of strictly negative sectional curvature? Does the metric converge to a hyperbolic one without forming singularities? 
 
 \begin{rem}\rm

 See Aitchison and Rubinstein \cite{AR1} for background information on
 cubings of non-positive curvature. Mosher \cite{Mo} proved that
 atoroidal 3--manifolds with such cubings have word hyperbolic
 fundamental groups. M Simon \cite{Si1} has studied the behaviour of
 singular metrics under Ricci flow. One might hope to find a similar
 polyhedral metric of non-positive curvature on atoroidal Haken
 3--manifolds, as in problem 16 below. If this could be done, one
 might hope to reprove Thurston's geometrisation theorem, using Ricci
 flow.

\end{rem}

\end{prob}

\begin{prob}

Do all hyperbolic 3--manifolds $M$ have cubical resolutions? These are
compact cubical complexes $\cal X$ which have CAT(0) structures on
their universal covers and are homotopy equivalent to $M$.  See
Rubinstein \cite{Ru2} and Rubinstein and Sageev \cite{RS1} for
discussion on how cubical resolutions can be constructed.

\end{prob}

\begin{prob}

Haken 3--manifolds have very short hierarchies (Aitchison and
Rubinstein \cite{AR2}). Can these be used in the atoroidal case to get
a combinatorial deformation theory to obtain CAT(0) or combinatorial
hyperbolic metrics?

 \begin{rem}\rm
 
Note that the gluing up of very short hierarchies is extremely
simple. One can take a handlebody $H$ and divide $\partial H$ into two
collections of subsurfaces, coloured black/white. This colouring has
to satisfy the condition that every meridian disk for $H$ has at least
two black and two white arcs in its boundary. Then one glues together
the black regions in pairs and the white regions become closed
incompressible surfaces. Finally the 3--manifolds formed from such
handlebodies are glued together along the closed incompressible
surfaces to form $M$. It turns out that the characteristic variety has
a simple description and so it is easy to detect whether the manifold
is atoroidal (cf Rubinstein \cite{Ru2}).  There is also a fairly
straightforward argument that the fundamental group of $M$ is word
hyperbolic (cf Rubinstein \cite{Ru2} and Swarup \cite{Sw}).

\end{rem}

\end{prob}

 \section {Haken 4--manifolds}
 
\begin{prob}[B Foozwell]
  
Show that every pair of Haken 3--manifolds $M, M^\prime$ are Haken cobordant, ie, there is a Haken 4--manifold $W$ so that $\partial W = M \cup M^\prime$ and the inclusions $M, M^\prime \subset W $ induce injections on fundamental groups. 

 \begin{rem}\rm

As a special case, given a Haken 3-manifold $M$, construct a Haken
 4--manifold $W$ so that $\partial W = M$ and the inclusion $M \subset
 W$ is an injection on fundamental groups.  Haken 4--manifolds are
 studied by Foozwell in \cite{Fo} and are built from hierarchies of
 $\pi_1$--injective 3--manifolds. Note the recent examples of Ivansic,
 Ratcliffe and Tschantz \cite{IR} which are Haken 4--manifolds
 obtained by removing a knotted 2--torus from the 4--sphere!
  
\end{rem}

\end{prob}
 
\begin{prob} 
 
Show that if $W_1, W_2$ are Haken 4--manifolds and there is a homotopy equivalence between $W_1$ and $W_2$, then $W_1$ and $W_2$ are homeomorphic. 
 
 \begin{rem}\rm
 
Note that this can be viewed as a weak version of Waldhausen's
classical result in dimension 3 (cf \cite{Wa}). However here we need
to assume that \emph{both} manifolds are Haken, whereas for
Waldhausen, only one is Haken and the other is assumed merely
irreducible. Note also that Haken 4--manifolds have `large'
fundamental groups, so Freedman's surgery theory cannot be applied.

\end{rem}

\end{prob}
 
\begin{prob} 
 
Study $\pi_1$--injective 3--manifolds in hyperbolic 4--manifolds. Many
 examples of hyperbolic 4--manifolds come from Coxeter group
 constructions, so naturally contain interesting immersed totally
 geodesic 3--manifolds which are certainly $\pi_1$--injective. However
 the construction in Aitchison, Matsumoto and Rubinstein \cite{AMR}
 appears possible to apply to some hyperbolic 4--manifolds with cusps,
 to construct many immersed 3--manifolds with principal curvatures at
 most $+1$, which are then $\pi_1$--injective. Also one could try to
 study Freedman--Freedman twisting \cite{FF} in dimension~4, assuming
 that suitable $\pi_1$--injective spanning 3--manifolds can be found
 for 4--manifolds with $\pi_1$--injective boundary.

\end{prob}
 
 \section {Contact structures and Heegaard splittings}
 
\begin{prob} 
  
Find a topological characterisation of the class of 3--manifolds which
admit the existence of tight contact structures. Many results have
been obtained on this problem---see Etnyre and Honda \cite{EH} and
Honda \cite{Ho1,Ho2}. There is a connection to the existence of taut
foliations, following the argument of Eliashberg and Thurston \cite{ET},
where a taut foliation is deformed to a tight contact structure.
  
 \end{prob}
 
\begin{prob} 
 
Giroux \cite{Gi} has shown that contact structures correspond to
certain types of open book decompositions of 3--manifolds. There is an
interesting connection between open book decompositions and Heegaard
splitting theory. We say that a Heegaard splitting $(M, H_+, H_-)$ can
be \emph{flipped} if there is an isotopy of $M$ interchanging the two
handlebodies $H_+, H_-$. Is it true that the Heegaard splitting can be
flipped if and only if the Heegaard surface $S = \partial H_+ =
\partial H_-$ consists of two pages of an open book decomposition, so
that the binding lies on $S$? Another version of this question is to
ask if the Heegaard surface $S$ has an infinite automorphism group and
$M$ is atoroidal, then is it true that $S$ comes from an open book
decomposition?

 \begin{rem}\rm

Note that in general, a Heegaard splitting might be two pages of
several open book decompositions. Equivalently, there could be several
fibred knots lying on the Heegaard surface dividing it into two
pages. (Clearly reducible splittings can be constructed with this
property.) In this case, there will be many flipping diffeomorphisms
obtained by composition, which do not correspond in a simple way to
open book structures. So this makes the problem quite
challenging. Note also that this problem is closely related to the
question of whether the natural map between the group $G$ of isotopy
classes of diffeomorphisms of the Heegaard surface, which extend over
the two handlebodies and the group $G^\prime$ of isotopy classes of
diffeomorphisms of $M$ has a finite or infinite kernel. For more
information on diffeomorphisms of $M$ which preserve a Heegaard
surface $S$, called automorphisms of the Heegaard splitting, see
Johnson and Rubinstein \cite{JR}.

 \end{rem}

\end{prob}
 
\begin{prob} 
 
Find  classes of closed orientable hyperbolic 3--manifolds, which only have a single irreducible or strongly irreducible Heegaard splitting, up to isotopy. Are such 3--manifolds common or rare?

 \begin{rem}\rm
 
 One likely class is obtained by large Dehn filling on the figure 8
 knot space. It is easy to see that the figure 8 knot space has only
 one irreducible splitting, since the only almost normal surfaces are
 obtained by adding an unknotted tube parallel to an edge to the
 peripheral torus. By the technique in Moriah and Rubinstein
 \cite{MR}, if one does suitably large Dehn filling, then any low
 genus irreducible/strongly irreducible splitting will actually be a
 splitting of the figure 8 knot space. Hence it is the unique
 splitting. However it is not clear that this is true for all
 splittings.

 \end{rem}

 \end{prob}
 
 \section {Singular incompressible surfaces}
 
\begin{prob} [W Thurston]
  
Does every closed or finite volume complete hyperbolic 3--manifold
admit an immersed or embedded closed surface with all principal
curvatures at most one? It is easy to show that such surfaces are
$\pi_1$--injective. (Compare with problem 1.) Note in Aitchison,
Matsumoto and Rubinstein \cite{AMR}, it is shown that the figure 8
knot space has a huge number of such surfaces.

 \begin{rem}\rm
 
  Here is a wonderful sketch of Thurston giving a way of studying this problem. Suppose we take a random point and a random totally geodesic plane $\Pi$ through that point $x$ in a closed hyperbolic 3--manifold $M$. As one takes an expanding domain in $\Pi$ starting at $x$, the boundary curve $\Gamma$ will grow exponentially, but will be in the bounded size manifold $M$. Consequently
if one could arrange that points in $\Gamma$ can be paired up to be reasonably close together, then by a small enough bending of $\Pi$, a new plane $P$ with small enough principal curvature could be found and a new domain and boundary curve $C$, which could be glued to itself to form a smooth surface satisfying the required conditions. 
  
 \end{rem}

 \end{prob}
  
 \begin{prob}[W Thurston]
 
Does every closed hyperbolic 3--manifold admit a finite sheeted cover
by a surface bundle over a circle? Note in Aitchison and Rubinstein
\cite{AR3}, a number of sets of examples of this type are
constructed. The basic idea is to glue together a collection of
fiberings of a fundamental domain. Can this method be made more
systematic, to produce a suitable `normal surface theory' of immersed
foliations?

 \begin{rem}\rm
 
  In \cite{AR3}, the vertex link structure must be a {\em regular
  branched cover}, to ensure that the singular foliation glues up to
  an immersion at each vertex. So for example, for a cubing of
  non-positive curvature, one needs to suppose that every edge has
  degree a multiple of 4 and every vertex link is a regular branched
  cover of an octahedron, to deduce that the manifold has an immersed
  fibration which lifts to a surface bundle structure in a finite
  sheeted cover. See Rubinstein \cite{Ru2} for further discussion of
  this necessary condition.

  \end{rem}

 \end{prob}
  
 \begin{prob} 

Suppose that a 3--manifold $M$ has a triple handlebody decomposition satisfying a disk condition. Does $M$ have an immersed incompressible surface without triple points? Are such surfaces separable, ie, do they lift to embeddings in finite sheeted coverings?

 \begin{rem}\rm
 
Triple handlebody decompositions are obtained by gluing three
handlebodies along subsurfaces of their boundaries. Triple curves are
then the boundaries between the different boundary subsurfaces. A disk
condition is the requirement that all meridian disks in the three
handlebodies intersect the triple curves at least $(m,n,p)$ times,
where ${1 \over m} +{1 \over n} +{1 \over p} \le {1 \over 2} $. In
\cite{Co}, numerous properties and constructions of 3--manifolds of
this type are given. The class can be viewed as generalising Seifert
fibred spaces with infinite fundamental groups and Haken
3--manifolds. One can view this question as the analogue of the
well-known fact that Seifert fibred spaces with infinite fundamental
groups all have immersed incompressible tori without triple points.

  \end{rem}

 \end{prob}

{\bf Acknowledgement}\qua The author is supported in part by the
Australian Research Council.

\bibliographystyle{gtart}
\bibliography{link}

\end{document}